\documentclass[12pt]{amsart}
\usepackage[all]{xy}

\newtheorem{thm}{Theorem}[section]
\newtheorem{prob}[thm]{Problem}

\newcommand{\sr}[2]{{\txt{$#1$\\$#2$}}}
\newcommand{\fx}{\mathfrak{x}}
\renewcommand{\b}{{\mathfrak b}}

\renewcommand{\d}{{\mathfrak d}}
\newcommand{\p}{{\mathfrak p}}
\newcommand{\s}{\mathfrak{s}}

\newcommand{\dnannouncement}[1]{[\S\ref{#1} below]}

\newcommand{\M}{\mathcal{M}}

\newcommand{\cov}{\mathsf{cov}}

\renewcommand{\b}{\mathfrak{b}}
\renewcommand{\t}{\mathfrak{t}}

\newcommand{\bq}{\begin{quote}}
\newcommand{\eq}{\end{quote}}
\renewcommand{\O}{\mathcal{O}}

\newcommand{\sone}{\mathsf{S}_1}    \newcommand{\sfin}{\mathsf{S}_{fin}}
\newcommand{\ufin}{\mathsf{U}_{fin}}

\newcommand{\nin}{\not\in}

\newcommand{\cU}{\mathcal{U}}

\newcommand{\sbst}{\subseteq}
\newcommand{\by}[2]{\par\hfill\emph{#1}, #2}
\newcommand{\Tau}{\mathrm{T}}
\newcommand{\CE}{\textsc{CE}}

\newcommand{\be}{\begin{enumerate}}
\newcommand{\ee}{\end{enumerate}}
\newcommand{\bi}{\begin{itemize}}
\newcommand{\ei}{\end{itemize}}
\renewcommand{\i}{\item}
\newcommand{\SPMBul}{\textbf{$\mathcal{SPM}$ Bulletin}}

\newcommand{\arx}[1]{\texttt{http://arxiv.org/abs/#1}}
\newcommand{\probmonth}{\emph{Problem of the month}}

\title[$\mathcal{SPM}$ Bulletin \textbf{4} (April 2003)]{%
$\mathcal{SPM}$ Bulletin\\[0.5cm]
Issue number 4: April 2003 \CE{}}

\begin{document}
\maketitle

\tableofcontents

\section{Editor's note}
This issue of the \SPMBul{} announces two conferences which are of interest to
anyone working in SPM or general topology. In the conference announced in
\dnannouncement{top} it is planned to have
\emph{a significant part devoted to SPM}.
Those who are interested in participating should contact
Ljubi\v{s}a D.\ R.\ Ko\v{c}inac at lkocinac@ptt.yu

Ko\v{c}inac is a very active mathematician in the field of SPM.
We announce here one of his most recent works.

The first issues of this bulletin,
which contain general information (first issue),
basic definitions, research announcements, and open problems (all issues) are available online:
\be
\i First issue: \arx{math.GN/0301011}
\i Second issue: \arx{math.GN/0302062}
\i Third issue: \arx{math.GN/0303057}
\ee
We are looking forward to receive more announcements from
the recipients of this bulletin.

\by{Boaz Tsaban}{tsaban@math.huji.ac.il}

\hfill \texttt{http://www.cs.biu.ac.il/\~{}tsaban}

\section{Research announcements}

\subsection{Selection principles in uniform spaces}
We begin the investigation of selection principles in uniform spaces in
a manner as it was done with selection principles theory for topological
spaces. We introduced and characterized uniform versions of classical
topological notions of the Menger, Hurewicz and Rothberger properties.
The uniform $\gamma$-sets are also considered.
\by{Ljubi\v{s}a D.\ R.\ Ko\v{c}inac}{lkocinac@ptt.yu}

\section{Conferences}

\subsection{General and Set-Theoretic Topology}
(dedicated to the 60th birthday of Istvan Juhasz)
Budapest, Hungary, from 8 to 13 August, 2003.

The J\'anos Bolyai Mathematical Society, the Alfr\'ed Renyi Institute of
Mathematics and the Paul Erd\~os Summer Research Center of Mathematics
are  organizing a Colloquium on General and Set-Theoretic Topology
in the period of August 8-13, 2003 in Budapest.  You are  cordially invited
to attend this conference. The aim of the Colloquium is to
provide ground for the exchange of information on new achievements and on the
recent problems of General and Set-Theoretic  Topology.

Organizing Committee:
    A.\ Csaszar (cochairman), J.\ Gerlits, A.\ Hajnal (chairman), E.\ Makai,
    G.\ Sagi, L.\ Soukup(secretary), Z.\ Szentmiklossy

The following people will give invited addresses:
\bi
\i A.\ V.\ Archangelskii  (Athens, OH)
\i A.\ Dow   (Charlton, NC)
\i K.\ Kunen      (Madison)
\i J.\ van Mill (Amsterdam)
\i S.\ Romaguera (Valencia)
\i S.\ Shelah (Jerusalem)
\i W.\ A.\ R.\ Weiss (Toronto)
\ei

Besides these lectures there will be 20-minutes contributed
talks.

The conference will take place at the  Alfred Renyi  Institute of
Mathematics, Budapest (V.\ Realtanoda u.\ 13--15, Budapest).

\hfill{\emph{E-mail:} top2003@renyi.hu

\hfill{\emph{URL:} \texttt{http://www.renyi.hu/\~{}top2003}

\subsection{Interplay between Topology and Analysis}\label{top}
(accompanying the Congress of MASSEE)
September 15--21, 2003, Bulgaria.
Organizing and Programme Committee:
\bi
\i A.\ V.\ Arhangel'skii, Univrsity of Ohaio, USA
\i M.\ M.\ Choban, Tiraspol University in Kishineu, Moldova
\i D.\ Dimovski, Univrsity of Skopje, Macedonia
\i P.\ S.\ Kenderov (Co-Chair), Institute of Mathematics and Informatics,  Bulgarian Academy of Sciences
\i Lj.\ Kocinac, University of Nis, Serbia
\i S.\ J.\ Nedev (Co-Chair), Institute of Mathematics and Informatics,  Bulgarian Academy of Sciences
\i J.\ Orihuela, University of Murcia, Spain
\i J.\ Pelant, Institute of Mathematics, Czech Academy of Sciences
\i I.\ Shishkov (Secretary), Institute of Mathematics and Informatics,  Bulgarian Academy of Sciences
\i St.\ Troyanski, Institute of Mathematics and Informatics,  Bulgarian Academy of Sciences
\ei

The main purpose of this Mini Symposium  is to bring together, within the Congress of
MASSEE, people interested in the indicated areas, including Ph.D.\ and university students.
Topology and the Analisys occupy central place in Mathematics. The interplay between Topology
and Analysis  is evident in such topics like Selections of set-valued mappings, Geometry of Banach
spaces and, in particular, equivalent renorming of Banach spaces.
Recent results in the area show that the renorming of Banach spaces is closely related to
s-fragmentability of the weak topologies, and thus to the topological theory of metrizability.
Useful tools when attacking problems in this area are Topological games and quasi-continuity.
Any other related topics are welcome.

The program of the Workshop will include invited talks of approximately 45 min.\ duration,
as well as contributed presentations of approximately 20 min.\ duration (including 5 min.\ discussion),
depending on the total number of contributions.
For the accommodation, participation fees, deadlines, abstracts, etc., please refer to
\bq
\texttt{http://www.math.bas.bg/massee2003}
\eq
\by{Ivailo Shishkov}{shishkov@math.bas.bg}

\section{Modern selection principles: Open problems}
By \emph{modern selection principles} we mean selection
principles involving modern types of open covers.
These include groupable and weakly groupable covers
\cite{coc7, coc8}, and $\tau$-covers.
In this note we concentrate on this latter type of covers.
The following is based on \cite{tautau}.

Fix a zero-dimensional, separable metric space $X$.
(Topologically this is simply a zero-dimensional set of reals.)
Recall that the symbols
$\O$, $\Omega$, and $\Gamma$ denote the by-now-classical
collections of (countable open)
covers, $\omega$-covers, and $\gamma$-covers of $X$.
$\cU$ is a \emph{$\tau$-cover} of $X$ if it is a large cover of $X$,
and for each $x,y\in X$, either $\{U\in\cU : x\in U, y\nin U\}$ is finite, or
$\{U\in\cU : y\in U, x\nin U\}$ is finite.
As any $\omega$-cover of a set of reals contains a countable $\omega$-cover
of that set, and every large subcover of a $\tau$-cover is a $\tau$-cover
(and every $\omega$-cover is a large cover), we have that
Every open $\tau$-cover of $X$ contains a countable $\tau$-cover of $X$.
Let $\Tau$ denote the collection of countable open $\tau$-covers of $X$.
Then
$$\Gamma  \sbst \Tau \sbst \Omega \sbst \O.$$
The notion of $\tau$-covers introduces seven new pairs---namely,
$(\Tau,\O)$, $(\Tau,\Omega)$, $(\Tau,\Tau)$, $(\Tau,\Gamma)$,
$(\O,\Tau)$, $(\Omega,\Tau)$, and $(\Gamma,\Tau)$---to
which any of the selection operators $\sone$, $\sfin$, and $\ufin$ can be applied.
This makes a total of $21$ new selection hypotheses.
Fortunately, some of them are easily eliminated,
and the surviving properties appear in the following diagram.

\medskip

{\scriptsize
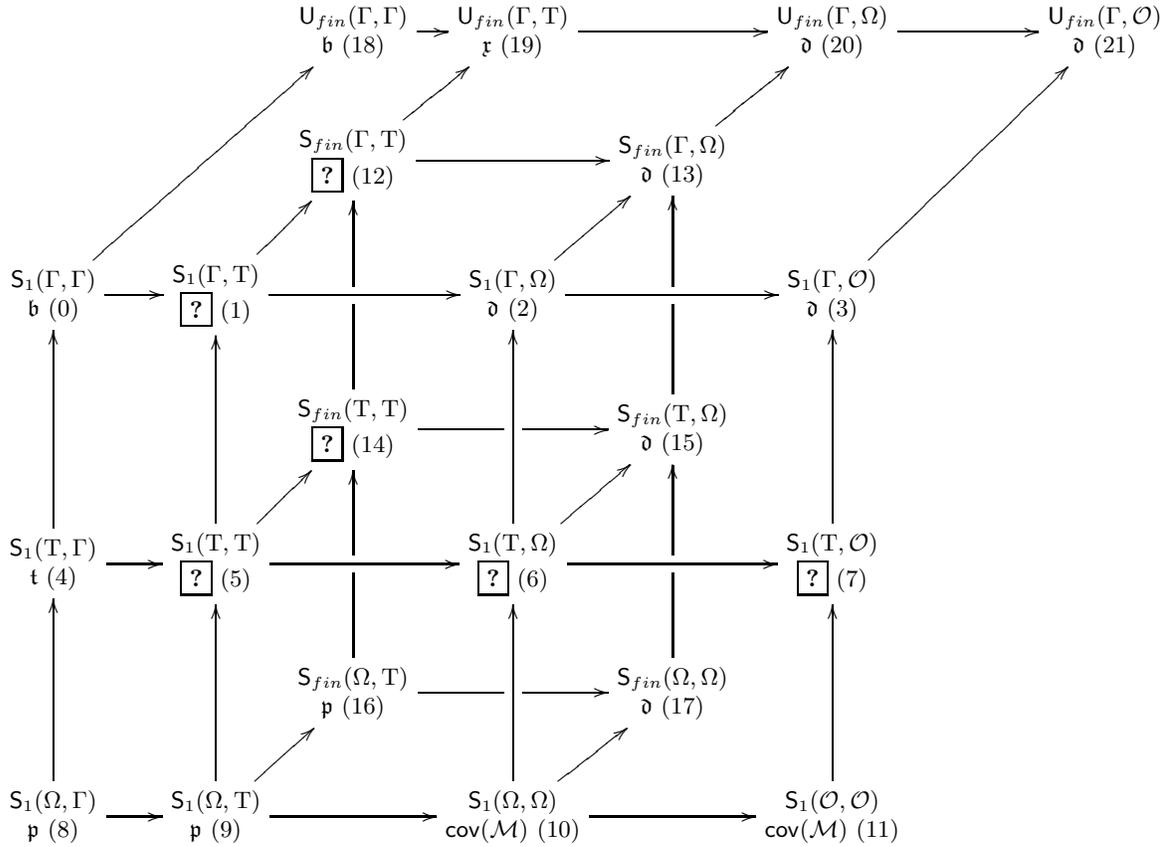
\begin{center}
$\xymatrix@C=7pt@R=20pt{
&
&
& \sr{\ufin(\Gamma,\Gamma)}{\b~~ (18)}\ar[r]
& \sr{\ufin(\Gamma,\Tau)}{\fx~~ (19)}\ar[rr]
&
& \sr{\ufin(\Gamma,\Omega)}{\d~~ (20)}\ar[rrrr]
&
&
&
& \sr{\ufin(\Gamma,\O)}{\d~~ (21)}
\\
&
&
& \sr{\sfin(\Gamma,\Tau)}{\fbox{\textbf{?}}~~ (12)}\ar[rr]\ar[ur]
&
& \sr{\sfin(\Gamma,\Omega)}{\d~~ (13)}\ar[ur]
\\
\sr{\sone(\Gamma,\Gamma)}{\b~~ (0)}\ar[uurrr]\ar[rr]
&
& \sr{\sone(\Gamma,\Tau)}{\fbox{\textbf{?}}~~ (1)}\ar[ur]\ar[rr]
&
& \sr{\sone(\Gamma,\Omega)}{\d~~ (2)}\ar[ur]\ar[rr]
&
& \sr{\sone(\Gamma,\O)}{\d~~ (3)}\ar[uurrrr]
\\
&
&
& \sr{\sfin(\Tau,\Tau)}{\fbox{\textbf{?}}~~ (14)}\ar'[r][rr]\ar'[u][uu]
&
& \sr{\sfin(\Tau,\Omega)}{\d~~ (15)}\ar'[u][uu]
\\
\sr{\sone(\Tau,\Gamma)}{\t~~ (4)}\ar[rr]\ar[uu]
&
& \sr{\sone(\Tau,\Tau)}{\fbox{\textbf{?}}~~ (5)}\ar[uu]\ar[ur]\ar[rr]
&
& \sr{\sone(\Tau,\Omega)}{\fbox{\textbf{?}}~~ (6)}\ar[uu]\ar[ur]\ar[rr]
&
& \sr{\sone(\Tau,\O)}{\fbox{\textbf{?}}~~ (7)}\ar[uu]
\\
&
&
& \sr{\sfin(\Omega,\Tau)}{\p~~ (16)}\ar'[u][uu]\ar'[r][rr]
&
& \sr{\sfin(\Omega,\Omega)}{\d~~ (17)}\ar'[u][uu]
\\
\sr{\sone(\Omega,\Gamma)}{\p~~ (8)}\ar[uu]\ar[rr]
&
& \sr{\sone(\Omega,\Tau)}{\p~~ (9)}\ar[uu]\ar[ur]\ar[rr]
&
& \sr{\sone(\Omega,\Omega)}{\cov(\M)~~ (10)}\ar[uu]\ar[ur]\ar[rr]
&
& \sr{\sone(\O,\O)}{\cov(\M)~~ (11)}\ar[uu]
}$
\end{center}
}

\medskip

Below each property in the above diagram appears a ``serial number'' (to be used below), and
its \emph{critical cardinality}, the minimal
cardinality of a set of reals not satisfying the property.
The cardinal numbers $\p$, $\t$, $\b$, and $\d$
are the well-known psedudo-intersection, tower, (un)bounding, and dominating numbers,
and $\cov(\M)$ is the covering number for the meager ideal
(see, e.g., \cite{vD} or \cite{BlassHBK} for definitions and details).
$\fx$ is the excluded-middle number, and is equal to $\max\{\s,\b\}$,
where $\s$ is the splitting number \cite{taucards}.

As indicated in the diagram, some of the critical cardinalities were not found yet.
\begin{prob}\label{cards}
What are the unknown critical cardinalities in the above diagram?
\end{prob}
Recall that Issue 2's \probmonth{} mentioned two unsettled implications in
the corresponding diagram for the \emph{classical} types of open covers.
As there are many more properties when $\tau$-covers are incorporated into the
framework, and since this investigation is new, there remain \emph{many} unsettled
implications in the above diagram. To be precise, there are exactly $76$(!) unsettled
implications in this diagram. These appear as question marks in the
following \emph{implications table}.
Entry $(i,j)$ in the table ($i$th row, $j$th column) is to be interpreted as follows:
It is $1$ if property $i$ implies property $j$, $0$ if property $i$ does not
imply property $j$ (that is, consistently there exists a counter-example),
and $?$ if the implication is unsettled.

\newcommand{\mb}[1]{{\mbox{\textbf{#1}}}}

\medskip

\begin{center}
{\scriptsize
\begin{tabular}{|r||cccccccccccccccccccccc|}
\hline
   & \mb{0} & \mb{1} & \mb{2} & \mb{3} & \mb{4} & \mb{5} & \mb{6} & \mb{7} & \mb{8} & \mb{9} & \mb{10} & \mb{11} & \mb{12} & \mb{13} & \mb{14} & \mb{15} & \mb{16} & \mb{17} & \mb{18} & \mb{19} & \mb{20} & \mb{21}\cr
\hline\hline
\mb{ 0} &  1 & 1 & 1 & 1 & 0 & \mb{?} & \mb{?} & \mb{?} & 0 & 0 & 0 & 0 & 1 & 1 & \mb{?} & \mb{?} & 0 & 0 & 1 & 1 & 1 & 1\cr
\mb{ 1} &  \mb{?} & 1 & 1 & 1 & 0 & \mb{?} & \mb{?} & \mb{?} & 0 & 0 & 0 & 0 & 1 & 1 & \mb{?} & \mb{?} & 0 & 0 & \mb{?} & 1 & 1 & 1 \cr
\mb{ 2} &  0 & 0 & 1 & 1 & 0 & 0 & \mb{?} & \mb{?} & 0 & 0 & 0 & 0 & 0 & 1 & 0 & \mb{?} & 0 & 0 & 0 & 0 & 1 & 1 \cr
\mb{ 3} &  0 & 0 & 0 & 1 & 0 & 0 & 0 & \mb{?} & 0 & 0 & 0 & 0 & 0 & 0 & 0 & 0 & 0 & 0 & 0 & 0 & 0 & 1 \cr
\mb{ 4} &  1 & 1 & 1 & 1 & 1 & 1 & 1 & 1 & 0 & 0 & \mb{?} & \mb{?} & 1 & 1 & 1 & 1 & 0 & \mb{?} & 1 & 1 & 1 & 1 \cr
\mb{ 5} &  \mb{?} & 1 & 1 & 1 & \mb{?} & 1 & 1 & 1 & 0 & 0 & \mb{?} & \mb{?} & 1 & 1 & 1 & 1 & 0 & \mb{?} & \mb{?} & 1 & 1 & 1 \cr
\mb{ 6} &  0 & 0 & 1 & 1 & 0 & 0 & 1 & 1 & 0 & 0 & \mb{?} & \mb{?} & 0 & 1 & 0 & 1 & 0 & \mb{?} & 0 & 0 & 1 & 1 \cr
\mb{ 7} &  0 & 0 & 0 & 1 & 0 & 0 & 0 & 1 & 0 & 0 & 0 & \mb{?} & 0 & 0 & 0 & 0 & 0 & 0 & 0 & 0 & 0 & 1 \cr
\mb{ 8} &  1 & 1 & 1 & 1 & 1 & 1 & 1 & 1 & 1 & 1 & 1 & 1 & 1 & 1 & 1 & 1 & 1 & 1 & 1 & 1 & 1 & 1 \cr
\mb{ 9} &  \mb{?} & 1 & 1 & 1 & \mb{?} & 1 & 1 & 1 & \mb{?} & 1 & 1 & 1 & 1 & 1 & 1 & 1 & 1 & 1 & \mb{?} & 1 & 1 & 1 \cr
\mb{10} &  0 & 0 & 1 & 1 & 0 & 0 & 1 & 1 & 0 & 0 & 1 & 1 & 0 & 1 & 0 & 1 & 0 & 1 & 0 & 0 & 1 & 1 \cr
\mb{11} &  0 & 0 & 0 & 1 & 0 & 0 & 0 & 1 & 0 & 0 & 0 & 1 & 0 & 0 & 0 & 0 & 0 & 0 & 0 & 0 & 0 & 1 \cr
\mb{12} &  \mb{?} & \mb{?} & \mb{?} & \mb{?} & 0 & \mb{?} & \mb{?} & \mb{?} & 0 & 0 & 0 & 0 & 1 & 1 & \mb{?} & \mb{?} & 0 & 0 & \mb{?} & 1 & 1 & 1 \cr
\mb{13} &  0 & 0 & 0 & 0 & 0 & 0 & 0 & 0 & 0 & 0 & 0 & 0 & 0 & 1 & 0 & \mb{?} & 0 & 0 & 0 & 0 & 1 & 1 \cr
\mb{14} &  \mb{?} & \mb{?} & \mb{?} & \mb{?} & \mb{?} & \mb{?} & \mb{?} & \mb{?} & 0 & 0 & \mb{?} & \mb{?} & 1 & 1 & 1 & 1 & 0 & \mb{?} & \mb{?} & 1 & 1 & 1 \cr
\mb{15} &  0 & 0 & 0 & 0 & 0 & 0 & 0 & 0 & 0 & 0 & 0 & 0 & 0 & 1 & 0 & 1 & 0 & \mb{?} & 0 & 0 & 1 & 1 \cr
\mb{16} &  \mb{?} & \mb{?} & \mb{?} & \mb{?} & \mb{?} & \mb{?} & \mb{?} & \mb{?} & \mb{?} & \mb{?} & \mb{?} & \mb{?} & 1 & 1 & 1 & 1 & 1 & 1 & \mb{?} & 1 & 1 & 1 \cr
\mb{17} &  0 & 0 & 0 & 0 & 0 & 0 & 0 & 0 & 0 & 0 & 0 & 0 & 0 & 1 & 0 & 1 & 0 & 1 & 0 & 0 & 1 & 1 \cr
\mb{18} &  0 & 0 & 0 & 0 & 0 & 0 & 0 & 0 & 0 & 0 & 0 & 0 & 0 & \mb{?} & 0 & \mb{?} & 0 & 0 & 1 & 1 & 1 & 1 \cr
\mb{19} &  0 & 0 & 0 & 0 & 0 & 0 & 0 & 0 & 0 & 0 & 0 & 0 & 0 & \mb{?} & 0 & \mb{?} & 0 & 0 & 0 & 1 & 1 & 1 \cr
\mb{20} &  0 & 0 & 0 & 0 & 0 & 0 & 0 & 0 & 0 & 0 & 0 & 0 & 0 & \mb{?} & 0 & \mb{?} & 0 & 0 & 0 & 0 & 1 & 1 \cr
\mb{21} &  0 & 0 & 0 & 0 & 0 & 0 & 0 & 0 & 0 & 0 & 0 & 0 & 0 & 0 & 0 & 0 & 0 & 0 & 0 & 0 & 0 & 1 \cr
\hline
\end{tabular}
}
\end{center}

\medskip

\begin{prob}
Settle any of the unsettled implications in the above table.
\end{prob}
Almost any solution logically implies several other solutions.
The above implications table, as well as the remarks which follow, were obtained
by using a simple computer program.

For example, if the solution to Issue 2's \probmonth{} is negative,
that is, $\ufin(\Gamma,\Gamma)$ does not imply $\sfin(\Gamma,\Omega)$
(a plausible assumption), then $6$ implications are settled.
If, in addition, the solution to Issue 1's \probmonth{} is
positive, then only $53$ implications (out of the $76$ we began with) remain unsettled.

Marion Scheepers asked us which single
solution would imply as many other solutions
as possible. The answer found by our program is the following:
If entry $(12,5)$ is $1$ (that is, $\sfin(\Gamma,\Tau)$ implies $\sone(\Tau,\Tau)$), then
there remain only 33 (!) open problems.
The best possible negative entry in $(16,3)$: If $\sfin(\Omega,\Tau)$ does not imply
$\sone(\Gamma,\O)$, then only $47$ implications remain unsettled.

Finally, observe that any solution in Problem \ref{cards}
would imply several new nonimplications.

\by{Boaz Tsaban}{tsaban@math.huji.ac.il}

\section{Problem of the month}
Scheepers chose the following problem out of all the problems discussed above
as the most interesting.
\begin{prob}
Does $\sone(\Omega,\Tau)$ imply the Hurewicz property $\ufin(\Gamma,\Gamma)$?
\end{prob}
The reason for this choice is that
if the answer is positive, then $\sone(\Omega,\Tau)$ implies
the Gerlitz-Nagy $(*)$ property \cite{GN}, which is equivalent to
another modern selection property.

Observe that a positive answer to Issue 1's \probmonth{} implies
a positive answer to this problem too.

\by{Boaz Tsaban}{tsaban@math.huji.ac.il}

\end{document}